\documentclass[11pt]{article}
\usepackage{amsmath, amsthm, amsfonts}
\usepackage{amssymb}
\usepackage{amscd}
\usepackage{verbatim}

\begin{document}
\newcommand{\half}{\frac12}
\newcommand{\eqdef}{\stackrel{{\mathrm def}}{=}}
\newcommand{\C}{{\mathbb C}}
\newcommand{\M}{{\mathcal M}}
\newcommand{\loc}{{\mathrm{loc}}}
\newcommand{\dx}{\,\mathrm{d}x}
\newcommand{\core}{C_0^{\infty}(\Omega)}
\newcommand{\sob}{W^{1,p}(\Omega)}
\newcommand{\sobloc}{W^{1,p}_{\mathrm{loc}}(\Omega)}
\newcommand{\merhav}{{\mathcal D}^{1,p}}
\newcommand{\be}{\begin{equation}}
\newcommand{\ee}{\end{equation}}
\newcommand{\mysection}[1]{\section{#1}\setcounter{equation}{0}}
%%%%%%%%%%%%%%%
\newcommand{\bea}{\begin{eqnarray}}
\newcommand{\eea}{\end{eqnarray}}
\newcommand{\bean}{\begin{eqnarray*}}
\newcommand{\eean}{\end{eqnarray*}}
\newcommand{\thkl}{\rule[-.5mm]{.3mm}{3mm}}
%%%%%%%%%%%%%%%%%%%%%%%%%%%
\newcommand{\cw}{\stackrel{\rightharpoonup}{\rightharpoonup}}
\newcommand{\id}{\operatorname{id}}
\newcommand{\supp}{\operatorname{supp}}
\newcommand{\wlim}{\mbox{ w-lim }}
\newcommand{\mymu}{{x_N^{-p_*}}}
\newcommand{\R}{{\mathbb R}}
\newcommand{\N}{{\mathbb N}}
\newcommand{\Z}{{\mathbb Z}}
\newcommand{\Q}{{\mathbb Q}}
\newcommand{\abs}[1]{\lvert#1\rvert}
%%%%%%%%%%%
\newtheorem{theorem}{Theorem}[section]
\newtheorem{corollary}[theorem]{Corollary}
\newtheorem{lemma}[theorem]{Lemma}
\newtheorem{definition}[theorem]{Definition}
\newtheorem{remark}[theorem]{Remark}
\newtheorem{proposition}[theorem]{Proposition}
\newtheorem{assertion}[theorem]{Assertion}
\newtheorem{problem}[theorem]{Problem}
%%%%%%%%%%%%%%%%%%
\newtheorem{conjecture}[theorem]{Conjecture}
\newtheorem{question}[theorem]{Question}
\newtheorem{example}[theorem]{Example}
\newtheorem{Thm}[theorem]{Theorem}
\newtheorem{Lem}[theorem]{Lemma}
\newtheorem{Pro}[theorem]{Proposition}
\newtheorem{Def}[theorem]{Definition}
\newtheorem{Exa}[theorem]{Example}
\newtheorem{Exs}[theorem]{Examples}
\newtheorem{Rems}[theorem]{Remarks}
\newtheorem{Rem}[theorem]{Remark}

\newtheorem{Cor}[theorem]{Corollary}
\newtheorem{Conj}[theorem]{Conjecture}
\newtheorem{Prob}[theorem]{Problem}
\newtheorem{Ques}[theorem]{Question}
\newcommand{\pf}{\noindent \mbox{{\bf Proof}: }}
\newcommand{\Lag}{{\mathcal L}}
%%%%%%%%%%%%%%%%%%
%\newenvironment{proof}{{\bf Proof.}}{\hfill $\bowtie$\vskip4mm}

\renewcommand{\theequation}{\thesection.\arabic{equation}}
\catcode`@=11 \@addtoreset{equation}{section} \catcode`@=12

%\begin{titlepage}

\title{Hardy inequalities for weighted Dirac operator}
\author{Adimurthi%\thanks{Research supported by a travel grant from the Swedish Research Council}
\\
 {\small Centre of Applicable Mathematics}\\ {\small
 Tata Institute of Fundamental Research}\\
 {\small P.O.Box No. 1234}\\
{\small Bangalore - 560 012, India}\\
{\small aditi@aditi@math.tifrbng.res.in}\\\and Kyril Tintarev
\\{\small Department of Mathematics}\\{\small Uppsala University}\\
{\small P.O.Box 480}\\
{\small SE-751 06 Uppsala, Sweden}\\{\small
kyril.tintarev@math.uu.se}}
%\date{}
 \maketitle
\newcommand{\dnorm}[1]{\thkl #1 \thkl\,}

\begin{abstract}
An inequality of Hardy type is established for quadratic forms involving Dirac operator and a weight $r^{-b}$ for functions in $\R^n$. The exact Hardy constant $c_b=c_b(n)$ is found and generalized minimizers are given. The constant $c_b$ vanishes on a countable set of $b$, which extends the known case $n=2$, $b=0$ which corresponds to the trivial Hardy inequality in $\R^2$.
Analogous inequalities are proved in the case $c_b=0$ under constraints and, with error terms, for a bounded domain.  
\\[2mm]
\noindent  2000  \! {\em Mathematics  Subject  Classification.}
Primary  \! 35Q40, 35Q75, 46N50, 81Q10; Secondary  35P05, 47A05, 47F05.\\[1mm]
\end{abstract}
%\end{titlepage}
%%%%%%%%%%%%%%%%%%%%%%%%%%%%%%%%%%%%%%%%%%%%%%%%%

\mysection{Introduction}
The well-known Hardy inequality in $\R^3$
$$
\int_{\R^3}|\nabla u|^2dx\ge\frac14\int_{\R^3}\frac{|u|^2}{|x|^2}dx
$$
expresses, in the context of classical quantum mechanics the celebrated uncertainty principle.
Since the seminal paper of Brezis and Vazquez \cite{BV}, Hardy inequality has received a renewed attention of numerous authors (here we quote only papers \cite{AA, AE, ACR, FT, ky2} that have an immediate connection to our results).

 Similar inequalities are known for relativistic versions of the Schr\"odinger equation, in particular, for the quadratic form of $\sqrt{-\Delta}$ (Kato inequality, \cite{Kato}), and for a quadratic form of Dirac operator $\int W(x)|(\sigma\cdot\nabla)u|^2dx$ for Pauli matrices $\sigma_i$, $i=1,2,3$, by Dolbeault, Esteban, S\'er\'e \cite{Dolb} and  Dolbeault, Esteban, Loss and  Vega \cite{Esteban}. The class of weights in the latter work has a specific decay rate at infinity.  

In this paper we prove Hardy inequality, with exact constants, for $W(x)=|x|^{-b}$, $b\in\R$, generalizing two known cases, $b=0$ corresponding to the usual Hardy inequality, and $b=-1$ corresponding to inequality (4) in \cite{Esteban}. The method, based on reduction to the usual weighted Sobolev inequalities in one dimension, allows to obtain similar estimates for general $W$. 
A surprising phenomenon observed here is a ``quantization of certainty'' - there is a discrete set of values $b$ for which the exact Hardy constant becomes zero.
We have opted here to prove the inequalities for functions in $C_0^\infty(\R^n\setminus\{0\})$, since this choice allows to consider all real values of $b$ rather than $b<n-2$, as well as to extend the range of parameters in the important Caffarelli-Kohn-Nirenberg inequalities (see Appendices A and B). We leave it as an exercise to the reader to show that for $b<n-2$ our inequalities (as well as Caffarelli-Kohn-Nirenberg inequalities for applicable parameters) follow from correspondent inequalities on $C_0^\infty(\R^n\setminus\{0\})$ by an elementary approximation argument (multiplication of the function near the origin by a family of cut-off functions).
\vskip3mm
Let $n\ge 2$,  $\sigma_i$, $i=1,\dots,m$, be Hermitian $m\times m$-matrices satisfying
\begin{equation}
\label{Clifford}
\sigma_i\sigma_j+\sigma_j\sigma_i=2\delta_{ij},\quad i,j=1,\dots,m.
\end{equation}
Such matrices are found, in particular, for $m=2^{n/2}$ when $n$ is even, and for 
$m=2^{(n+1)/2}$ when $n$ is odd. In particular, for $n=3$ one usually fixes the set of $\sigma_i$
as Pauli matrices
$$
\sigma_1=\left(\begin{array}{cc}
0 & 1\\
1 & 0\\
\end{array}\right),\quad
\sigma_2=\left(\begin{array}{cc}
0 & -i\\
i & 0\\
\end{array}\right),\quad
\sigma_3=\left(\begin{array}{cc}
1 & 0\\
0 & -1\\
\end{array}\right).
$$

As we study here equations for functions $\R^n\to \C^m$, we will use distinct notations for the scalar products in the domain and in the range of the functions: $\langle f,g\rangle$ for the scalar product in $\C^m$ and
$p\cdot q$ for the scalar product in $\R^n$.
The weighted Dirac operator is induced by the following quadratic form on $C_0^\infty(\R^n\setminus\{0\};\C^m)$
\begin{equation} 
\label{Qb}
Q_b(u)\eqdef \int_{\R^n} r^{-b} |(\sigma\cdot\nabla) u|^2dx, b\in\R. 
\end{equation}
Quadratic form $Q_b(u)$ endows $C_0^\infty(\R^n\setminus\{0\};\C^m)$ with a scalar product, but the completion of $C_0^\infty(\R^n\setminus\{0\};\C^m)$ with respect to the corresponding norm is generally not a function space. 
An inequality of Hardy type, which is the main objective of this paper, yields an imbedding of the completed Hilbert space into a weighted $L^2$-space. The imbedding fails for a countable subset of $b$, in which case an analogous inequality holds on a subspace of $C_0^\infty(\R^n\setminus\{0\};\C^m)$.

The operator $\sigma\cdot\nabla$ is well known as a ``square root'' of the Laplacian, i.e.  
$(\sigma\cdot\nabla)^2=\Delta$, and it admits the following representation in polar coordinates (see e.g. \cite{Jerison}):
\begin{equation} 
\label{polar}
(\sigma\cdot\nabla) u= (\hat x\cdot\sigma)(\hat x\cdot\nabla u+\frac{1}{|x|}Lu),
\end{equation}
where $\hat x\eqdef \frac{x}{|x|}$, and the operator
\begin{equation} 
\label{L}
L\eqdef \sum_{j<k} \sigma_j\sigma_k(x_j\partial_{x_k}-x_k\partial_{x_j})
\end{equation}
involves differentiation only in the directions tangential to a sphere centered at the origin, that is, $[L,\partial_r]=0$. 
An elementary calculation based on evaluation of the integral of the squared magnitude in the right and the left hand side of \eqref{polar} shows that 
\begin{equation} 
\label{DeltaS}
L^2-(n-2)L=-\Delta_S,
\end{equation}
where $\Delta_S$ is the Laplace-Beltrami operator on the sphere. It is known (see e.g. \cite{Jerison}) that the spectrum $S_L$ of $L$ is discrete and consists of integer values:
\begin{equation} 
S_L=\{\Z\setminus\{1,\dots,n-2\}\}.
\end{equation}
We will denote the eigenspace of $L$ corresponding to the eigenfunction $k\in S_L$ as  $E_k$.

Our main results are as follows.

\begin{theorem}
\label{thm2-hardy}
Let $b\in\R$. Then for all $u\in C_0^\infty(\R^n\setminus\{0\};\C^m)$,
\begin{equation}
\label{Hardy}
\int_{\R^n} r^{-b} |(\sigma\cdot\nabla) u|^2dx- c_b\int_{\R^n} r^{-b-2} |u|^2dx\ge 0,
\end{equation}
where 
\begin{equation}
\label{cb}
 c_b\eqdef \min_{k\in\Z\setminus\{1,\dots,n-2\}}\left(k-\frac{n-2-b}{2}\right)^2
\end{equation}
is the exact constant for the inequality. Furthermore, the quadratic form in the left hand side
of \eqref{Hardy} has a finite-dimensional space of generalized ground states in the sense of Definition~\ref{def:gs}, that in polar coordinates have the form 
$v(r,\omega)=r^\frac{b+2-n}{2}\varphi(\omega)$, $\varphi\in E^{(b)}$, where $E^{(b)}$ is the span of (at most two) eigenspaces $E_k$ corresponding to those $k$ that yield the minimum in \eqref{cb}.  
\end{theorem}
Note that the maximal value of the constant $c_b$ is $\left(\frac{n-2}{2}\right)^2$, attained when $b=0$ and the inequality \eqref{Hardy} becomes the usual Hardy inequality. The constant $c_b$ equals zero if and only if $\frac{n-2-b}{2}\in\Z\setminus\{1,\dots,n-2\}$. This vanishing is not entirely unexpected, as this occurs also in the known case of Hardy inequality in two dimensions ($b=0$ and $n=2$).
In classical quantum mechanics, Hardy inequality (our case $b=0$) expresses the uncertainty principle (which remarkably fails for $n=2$).

\begin{theorem}
\label{thm:Sob}
Let $n>2$ and $2^*=2n/(n-2)$. Assume that $\frac{n-2-b}{2}\notin\Z\setminus\{0\dots,n-2\}$. Then 
there is a $C>0$ dependent on $n$ and $b$, such that for every $v\in C_0^\infty(\R^n\setminus\{0\})$,
\begin{equation}
\label{weightedSobolev1}
\int_{\R^n}|x|^{-b}|(\sigma\cdot\nabla) u|^2 dx \ge C\left(\int_{\R^n}|x|^{-\beta}|u|^{2^*} dx\right)^{2/2^*},
\end{equation}
where $\beta=\frac{b n}{n-2}$.
\end{theorem}
Note that the case $b=0$ gives the usual Sobolev inequality.

\mysection{Proof of Hardy inequality}
In this section we derive a representation of the quadratic form \eqref{Qb} in polar coordinates and prove Theorem~\ref{thm2-hardy}
\begin{lemma}
\label{representation-1}
For every $u\in C_0^\infty(\R^n\setminus\{0\};\C^m)$, 
\begin{equation}
\label{qba3}
Q_b(u)=\int_{\R^n}r^{-b}\left(|\partial_r u|^2+r^{-2}\langle (L^2+(b-n+2)L)u,u\rangle\right).
\end{equation}
\end{lemma}
\begin{proof} By \eqref{polar} and taking into account that $(\hat x\cdot\sigma)^2=I_m$ we have
\begin{equation}
\label{qba2}
Q_b(u)=\int_{\R^n}r^{-b}\left(|\partial_r u|^2+r^{-1}\partial_r\langle Lu,u\rangle + r^{-2}\langle L^2u,u\rangle\right).
\end{equation}
The second term in the right hand side can be evaluated by partial integration:
\begin{eqnarray*}
%\label{qb2a}
\int_{\R^n}r^{-b-1}\partial_r\langle Lu,u\rangle dx=\int_0^\infty\int_{S^{n-1}}r^{-b+n-2}\partial_r\langle Lu,u\rangle dr d\omega
\\= (b-n+2)\int_0^\infty\int_{S^{n-1}}r^{-b+n-3}\langle Lu,u\rangle drd\omega=(b-n+2)\int_{\R^n}r^{-b-2}\langle Lu,u\rangle dx.
\end{eqnarray*}
Substituting this into \eqref{qba2} we obtain \eqref{qba3}.
\end{proof}
\bigskip
{\em Proof of Theorem~\ref{thm2-hardy}}. Let us expand the function $u$ in the eigenfunctions of operator $L$, normalized in $L^2(S^{n-1})$:
$$
u=\sum_{k\in S_L}c_k(r)\psi_k(\omega).
$$
In order to account for the multiplicity of the eigenvalues we understand the term $c_k(r)\psi_k(\omega)$ as an implicit finite sum over a basis in the $k$-th eigenspace of $L$. 
We have from \eqref{qba3}
\begin{equation}
\label{qba4}
Q_b(u)=\sum_{k\in S_L}\int_0^\infty r^{-b+n-1}\left(|\partial_r c_k(r)|^2+r^{-2}(k^2+(b+2-n)k|c_k|^2\right)dr.
\end{equation}
Let us apply one-dimensional Caffarelli-Kohn-Nirenberg inequality \eqref{eq:CKN} to the first term in the right hand side. Note that the best constant in \eqref{eq:CKN} does not increase when one replaces $C_0^\infty(\R\setminus\{0\})$ with $C_0^\infty((0,\infty))$. 
\begin{equation}
\label{qba5}
\int_0^\infty r^{-b+n-1}|\partial_r c_k(r)|^2dr\ge \left(\frac{-b+n-2}{2}\right)^2\int_0^\infty r^{-b+n-3}|c_k(r)|^2dr.
\end{equation}
Substitution of \eqref{qba5} into \eqref{qba4} and collection of similar terms gives immediately  
\begin{equation}
\label{qba6}
Q_b(u)\ge\sum_{k\in S_L} \left(k+\frac{b+2-n}{2}\right)^2 \int_0^\infty r^{-b-2}|c_k|^2r^{n-1}dr.
\end{equation}
Inequality \eqref{Hardy} follows immediately.

It remains to show that the constant is exact. Indeed, let $k\in S_L$ be a value 
where the minimum of $\left(k-\frac{-b+n-2}{2}\right)^2$ is attained. Then, since we used the exact constant for \eqref{eq:CKN}, a minimizing sequence for \eqref{Hardy} is given by $\{c_j(r)\psi_k(\omega)\}_j$, where $\psi_k$ is an eigenfunction of $L$ with the eigenvalue  $k$, and  $c_j(r)$ is a minimizing sequence for \eqref{eq:CKN}, which with necessity converges  
to a scalar multiple of $r^\frac{b+2-n}{2}$ uniformly on compact subsets of $(0,\infty)$ (see \cite{ky2} for details.) 
\qed

\mysection{Proof of the Sobolev inequality}

{\em Proof of Theorem~\ref{thm:Sob}}
Let $\epsilon\in(0,1)$ be a constant to be specified at a later step.
From \eqref{qba3}, using \eqref{DeltaS}, we have
\begin{eqnarray*}
Q_b(u)=\epsilon \int_{\R^n}r^{-b}|\nabla u|^2dx
+\int_{\R^n}r^{-b}&\times&
\\
\left(
(1-\epsilon)|u_r|^2
+(1-\epsilon)r^{-2}\langle L^2u,u\rangle
%\\
+r^{-2}[(b-n+2)-\epsilon(n-2)]
\langle Lu,u\rangle
\right)dx.
\end{eqnarray*}
Expanding the expression under the second integral in the eigenfunctions of $L$ we obtain
\begin{eqnarray*}
Q_b(u)=\epsilon \int_{\R^n}r^{-b}|\nabla u|^2dx+\sum_{k\in S_L}\int_{\R^n}r^{-b}&\times&
\\
\left(
(1-\epsilon)|c_k'|^2+(1-\epsilon)r^{-2}k^2|c_k|^2
+r^{-2}[(b-n+2)-\epsilon(n-2)]
k|c_k|^2\right)
r^{n-1}dr.
\end{eqnarray*}
We apply the one-dimensional inequality \eqref{eq:CKN}:
$$
\int_{\R^n}r^{-b}|c_k'|^2r^{n-1}dr\ge \left(\frac{-b+n-2}{2}\right)^2\int_{\R^n}r^{-b-2}|c_k|^2r^{n-1}dr,
$$
and subsequently,
\begin{eqnarray*}
Q_b(u)=\epsilon \int_{\R^n}r^{-b}|\nabla u|^2dx
+\sum_{k\in S_L}\int_{\R^n}r^{-b-2}|c_k|^2&\times&
\\
(1-\epsilon)\left\{
\left(\frac{-b+n-2}{2}\right)^2
+k^2+[(b-n+2)-\epsilon(n-2)/(1-\epsilon)]
k\right\}
r^{n-1}dr.
\end{eqnarray*}
We can estimate the expression inside the large braces as follows, 
denoting $h(t)= \frac{t-n+2}{2}$ and $b_\epsilon\eqdef b-\epsilon(n-2)/(1-\epsilon)$:
\begin{eqnarray*}
h(b)^2+k^2+2h(b)k-\epsilon(n-2)/(1-\epsilon)k=
\\
=h(b)^2+k^2+2h(b_\epsilon)k\\
=(k+h(b_\epsilon))^2-h(b_\epsilon)^2+h(b)^2.
\ge C_{b_\epsilon}-C\epsilon.
\end{eqnarray*}
Since $c_b>0$ by assumption on $b$, and the last expression is continuous with respect to $\epsilon$, we choose fix an $\epsilon>0$ sufficiently small so that the expression is positive. This immediately implies that 
$$
Q_b(u)\ge\epsilon \int_{\R^n}r^{-b}|\nabla u|^2dx,
$$
and equation \ref{weightedSobolev1} follows immediately from \ref{weightedSobolev}, once we note that the condition on the exponent of $r$ for the latter inequality is $b\neq n-2$, which is included, as 
$\frac{n-2-b}{2}\neq 0$, into the assumptions of the theorem.
\qed

\mysection{Further inequalities}
Let $P_k$, $k\in S_L$, denote the orthogonal projectors in $L^2(\R^n;\C^m)$ induced by orthogonal projection on the eigenspace $E_k$ of $L$ in $L^2(S^{n-1};\C^m)$  
Note that while $\mathrm{dim}E_k<\infty$, $\mathrm{dim}P_kL^2(\R^n;\C^m)=+\infty$.

\begin{theorem}
\label{thm3-hardy'}
Assume that $j=\frac{n-2-b}{2}\in \Z\setminus\{1,\dots,n-2\}$. Then for all $u\in C_0^\infty(\R^n\setminus\{0\};C^m)$ satisfying $P_j u=0$, 
\begin{equation}
\label{Hardy-appended}
Q_b(u)=\int_{\R^n} r^{-b} |(\sigma\cdot\nabla) u|^2dx\ge \int_{\R^n} r^{-b-2} |u|^2dx,
\end{equation}
and the constant $1$ in the right hand side cannot be improved.
\end{theorem}
\begin{proof}
Consider relation \eqref{qba6} under the orthogonality condition $P_j u=0$, that is, assuming that the term corresponding to $k=j$ is excluded from the sum representing $Q_b(u)$ in \eqref{qba6}. Then $k=j-1$ or $k=j+1$ lies in $S_L$ and the smallest coefficient $(k-j)^2$ remaining in the expansion of \eqref{qba6} corresponds to $k=j\pm 1$ and equals $1$. An argument repetitive of that in the proof of Theorem~\ref{thm2-hardy} shows that this coefficient is exact. 
\end{proof}

An analog of the theorem above using a finite-dimensional projector gives a somewhat weaker inequality, connected the fact that $\int|\nabla u|^2$ cannot dominate any weighted $L^2$-norm in $\R^2$, and that in restriction to subspace of codimension 1 it still cannot dominate the $L^2$-norm with the Hardy weight $r^{-2}$. 

\begin{theorem}
\label{thm4-hardy'}
Assume that $j=\frac{n-2-b}{2}\in \Z\setminus\{1,\dots,n-2\}$. Then there is a $C>0$ such that for all $u\in C_0^\infty(\R^n\setminus\{0\};C^m)$ satisfying 
\begin{equation}
\int_1^2\int_{S^{n-1}} P_j u(r,\cdot)rdrd\omega=0,
\end{equation} 
the following inequality holds true:
\begin{equation}
\label{Hardy-appended4}
Q_b(u)=\int_{\R^n} r^{-b} |(\sigma\cdot\nabla) u|^2dx\ge C \int_{\R^n} r^{-b-2}(1+|\log r|)^{-2} |u|^2dx.
\end{equation}
\end{theorem}
\begin{proof}
Repeating the proof of Theorem~\ref{thm3-hardy'}, it suffices to estimate from below, 
under the orthogonality condition 
\begin{equation}
\int_1^2c_j(r)dr=0,
\end{equation} the term in the expansion \eqref{qba4}
corresponding to $k=j$, that is, the expression
\begin{equation}
\label{i1}
I\eqdef\int_0^\infty r^{-2j}\left(|\partial_r c_j(r)|^2-j^2|c_j|^2\right)rdr.
\end{equation}
Substituting $w(r)=c_j(r)r^j$ we get 
\begin{equation}
\label{i2}
I=\int_0^\infty |\partial_r w(r)|^2rdr.
\end{equation}
Then by \eqref{eq:H2D}
\begin{equation}
\label{i3}
I\ge C\int_0^\infty r^{-2}(1+|\log r|)^{-2}|w(r)|^2rdr=C\int_0^\infty r^{-b-2}(1+|\log r|)^{-2}|c_j|^2r^{n-1}dr.
\end{equation}
This yields \eqref{Hardy-appended4}.
\end{proof}

\begin{theorem}
\label{thm5-hardy'}
Let $\Omega\subset\R^n$ be a bounded domain containing the origin and let $R=\sup_{x\in\Omega}|x|$.
Let $\eta_1(r)\eqdef \log(R/r)$ and define recursively $\eta_j(r)\eqdef \eta_1\circ\eta_{j-1}$. 
Then for all $u\in C_0^\infty(\Omega\setminus\{0\};C^m)$,
\begin{eqnarray}
\label{Hardy-error}
\nonumber
Q_b(u)=\int_{\Omega} r^{-b} |(\sigma\cdot\nabla) u|^2dx\ge c_b\int_{\Omega} r^{-b-2}|u|^2dx 
\\
+c_b\sum_{k=1}^\infty.\int_{\Omega} r^{-2}\eta_1(r)...\eta_k(r) |u|^2dx 
\end{eqnarray}
\end{theorem}
\begin{proof}
The proof is repetitive of that of Theorem~\ref{thm2-hardy}, except for instead of substituting into \eqref{qba4} the estimate \eqref{qba5}, one uses a refinement of \eqref{qba5} for sunctions supported on a ball $B_R(0)$.
Modifying the reduction to the standard Hardy inequality by means of Theorem~\ref{thm:ckn}, one replaces \eqref{eq:Hardy} with the following well-known Hardy inequality with remainder terms (\cite{FT},\cite{ACR}):
\begin{equation}
 \label{eq:Hardy-err} \int_0^1 |w'|^2dr\ge \frac14 
\int_0^R r^{-2}|w|^2dr + \frac14\sum_{k=1}^\infty\int_0^Rr^{-2}\eta_1(r)\dots \eta_k(r)|w|^2dr
\end{equation}
\end{proof}
Several other error expressions refining \eqref{eq:Hardy} found in literature can be used to provide refinements of \eqref{qba5}, and subsequently, of \eqref{Hardy} as well. We leave it as an exercise for the reader.

\section*{Appendix A: Caffarelli-Kohn-Nirenberg inequality}
 
\begin{theorem}
\label{thm:ckn}
 Let $a\in\R$, $n\in\N$. Then for every $u\in C_0^\infty(\R^n)\setminus\{0\}$,
\begin{equation}
 \label{eq:CKN} \int_{\R^n} |x|^{a}|\nabla u|^2dx\ge \left(\frac{a+n-2}{2}\right)^2\int_{\R^n} |x|^{a-2}|u|^2dx,
\end{equation}
and the constant in the right hand side is exact. 
\end{theorem}
When $a>2-n$ the inequality is due to \cite{ckn}, where it is established for 
all $u\in C_0^\infty(\R^n)$. However once one considers only the functions vanishing near the origin the restriction  $a>2-n$ can be be removed and a significantly shorter proof can be given.  
\begin{proof}
The following well-known identity holds for any positive $\psi\in C^2(\R^n\setminus\{u\})$ and any $u\in C_0^\infty(\R^n)\setminus\{0\}$: 
\begin{equation}
 \label{eq:CKN1} \int_{\R^n} \psi^2 |\nabla u|^2dx=\int_{\R^n} (|\nabla(u\psi)|^2+(\Delta \psi)\psi|u|^2dx.
\end{equation}
Applying it to $\psi=|x|^{a/2}$, we have 
\begin{equation}
 \label{eq:CKN2} \int_{\R^n} |x|^{a}| |\nabla u|^2dx=
\int_{\R^n} \left( |\nabla(u\psi)|^2+\frac{a}{2}(a/2+n-2)|x|^{a-2}|u|^2\right)dx.
\end{equation}
Applying to the first term the standard Hardy inequality with the exact constant 
\begin{equation}
 \label{eq:Hardy} \int_{\R^n} |\nabla w|^2dx\ge \left(\frac{n-2}{2}\right)^2
\int_{\R^n}|x|^{-2}|w|^2dx
\end{equation}
and collecting similar terms in \eqref{eq:CKN2}, we arrive at \eqref{eq:CKN}.
\end{proof}
\section*{Appendix B: Weighted Sobolev inequality}
The inequality below for $\alpha>2-n$ is due to \cite{ckn} and \cite{Mazya}[2.1.6, Cor.2].
For a smaller class of functions $C_0^\infty(\R^n\setminus\{0\})$ the inequality also extends to $\alpha<2-n$.
\begin{proposition}
Let $n>2$ and $2^*=2n/(n-2)$. There is a $C>0$ such that for every $v\in C_0^\infty(\R^n\setminus\{0\})$,
\begin{equation}
\label{weightedSobolev}
\int_{\R^n}|x|^\alpha|\nabla v|^2 dx \ge C\left(\int_{\R^n}|x|^\beta|v|^{2^*} dx\right)^{2/2^*},
\end{equation}
where $\alpha\in\R\setminus\{2-n\}$  and $\beta=\frac{\alpha n}{n-2}$.
\end{proposition}
\begin{proof}
By Sobolev inequality, there exists $C>0$ such that for all $v\in C_0^\infty(\R^n\setminus\{0\})$,
\begin{equation}
\left(\int_{1<|x|\le 2}|v|^{2^*} dx\right)^{2/2^*}\le C \int_{1< |x|\le 2}|\nabla v|^2dx +
C \int_{1< |x|\le 2}|v|^2dx. 
\end{equation}
Then, with a renamed $C$, 
\begin{equation}
\left(\int_{1<|x|\le 2}|x|^\beta|v|^{2^*} dx\right)^{2/2^*}\le C \int_{1< |x|\le 2}|x|^\alpha|\nabla v|^2dx +
C \int_{1< |x|\le 2}|x|^{\alpha-2}|v|^2dx. 
\end{equation}
Replacing $v$ with $2^{\frac{j(N+\alpha-2)}{2}}v(2^j\cdot)$, $j\in\Z$, we obtain 
\begin{equation}
\left(\int_{2^j<|x|\le 2^{j+1}}|x|^\beta|v|^{2^*} dx\right)^{2/2^*}\le C \int_{2^j<|x|\le 2^{j+1}}|x|^\alpha|\nabla v|^2dx +
C \int_{2^j<|x|\le 2^{j+1}}|x|^{\alpha-2}|v|^2dx. 
\end{equation}
Addition over $j\in\Z$ provides, taking into account subadditivitiy in the left hand side, 
 \begin{equation}
\left(\int_{\R^n}|x|^\beta|v|^{2^*} dx\right)^{2/2^*}\le C \int_{\R^n}|x|^\alpha|\nabla v|^2dx +
C \int_{\R^n}|x|^{\alpha-2}|v|^2dx. 
\end{equation}
By \eqref{eq:CKN} the second term in the right hand side is dominated by the first term, from which, once we note that $\alpha\neq 2$, \eqref{weightedSobolev} is immediate.
\end{proof}
\section*{Appendix C: Hardy inequality in $\R^2$}
\begin{theorem}
\label{thm:H2D} 
Let $\Psi(u)=\int_{1<|x|<2}u(x)dx$. Then 
there exists $C>0$ such that for ever $u\in C_0^\infty(\R^2\setminus\{0\})$ satisfying $\Psi (u)=0$, the following inequality holds true: 
\begin{equation}
\label{eq:H2D}
\int_{\R^2} |\nabla u|^2 dx \ge C\int_{\R^2}r^{-2}(1+|\log r|)^{-2}|u|^2dx.
\end{equation}
\end{theorem}
\begin{proof} 
Assume first that $u\in C_0^\infty(B_R(0)\setminus\{0\})$, $R>0$ and is radially symmetric. Then it is easy to show that \eqref{eq:H2D} holds unconditionally by using the change of variable $t=\log(2R/|x|)$, which reduces the inequality to the standard one-dimensional Hardy inequality. This argument can be repeated for radially symmetric functions $u\in u\in C_0^\infty(\R^2\setminus B_1(0))$ (or by applying Kelvin tranformation to the inequality in the ball). From this and the elementary density argument follows that \eqref{eq:H2D} holds for all radially symmetric functions $u\in C_0^1(\R^2\setminus\{0\})$ satisfying $u(R)>0$ or $\Psi_0 (u)=0$, 
where $\Psi_0(u)=\int_{|x|=R}udx$. 
This can be equivalently rewritten (with a different constant $C$) as 
\begin{equation}
\label{H2D1}
\Psi_0(u)^2+\int_{\R^2} |\nabla u|^2 dx \ge C\int_{\R^2}r^{-2}(1+|\log r|)^{-2}|u|^2dx
\end{equation}
for all radially symmetric functions $u\in C_0^1(\R^2\setminus\{0\})$.
Assume now that
\begin{equation}
\label{H2D2}
\inf\{ \int_{\R^2} |\nabla u|^2:\quad  \int_{\R^2}r^{-2}(1+|\log r|)^{-2}|u|^2dx=1,\Psi(u)=0\}=0,
\end{equation}
where the infimum is taken over all radially symmetric functions $u\in C_0^1(\R^2\setminus\{0\})$ satisfying the constraints.
Then there exists a sequence of radial functions $u_k\in C_0^1(\R^2\setminus\{0\})$ such that 
$\int_{\R^2} |\nabla u_k|^2\to 0$, $\int_{\R^2}r^2(1+|\log r|)^{-2}|u_k|^2dx=1$ and $\Psi(u_k)=0$. Then $u_k$ is bounded in 
$H^1_{\mathrm loc}(\R^2)$ and, on a renumbered subsequence, it converges weakly in $H^1_{\mathrm loc}(\R^2)$ to some constant $\lambda$ satisfying $\Psi(\lambda)=0$. Then $\lambda=0$, which implies that  
$\Psi_0(u_k)\to 0$, and therefire it follows from \eqref{H2D1} that $\int_{\R^2}r^2(1+|\log r|)^{-2}|u_k|^2dx\to 0$, a contradiction. 
From \eqref{H2D2} follows 
\begin{equation}
\label{H2D3}
\Psi(u)^2+\int_{\R^2} |\nabla u|^2 dx \ge C\int_{\R^2}r^{-2}(1+|\log r|)^{-2}|u|^2dx
\end{equation}
for all radially symmetric functions $u\in C_0^1(\R^2\setminus\{0\})$.
The inequality for general functions $u\in C_0^1(\R^2\setminus\{0\})$ follows from decomposition into spherical harmonics once one observes that $\Psi(u)$ depends only the radial component of $u$ and that the the spherical Laplacian for any higher harmonics yields the lower bound of the form $C\int_{\R^2}r^{-2}|u|^2dx$.
\end{proof}
\section*{Appendix D: Definition of generalized ground state}
We give a generalization of definition of \cite{ky2} for real valued functionals on the normed (not necessarily complete) vector space. 
\begin{definition}
\label{def:gs}
 Let $X$ be a normed vector space and 
assume that a functional $Q:X\to \R$ is positively homogogeneous of a positive degree and that 
\begin{equation}
\label{gs} 
\inf_{\langle \xi,x\rangle=1} Q(x)=0
\end{equation}
for every $\xi\in X^*$. One says that $v\in X^{**}$
is a generalized ground state of $Q$ if for every $\xi\in X^*$ there exists a minimizing sequence $x_k$ for \eqref{gs}, whose weak-* limit is a scalar multiple of $v$. 
\end{definition}

This definition allows to define a ground state when \eqref{gs} has no minimizer even under an appropriate extension of $Q$ (the functional $Q$ is not required to be semicontinuous and $X$ is not required to be complete.)
When $Q$ is a positive quadratic form of the Schr\"odinger operator with a potential term, defined on $C_0^\infty(\Omega)$, the classical ground state of $Q$ in $L^2(\Omega)$ is also the ground state in the sense of the definition above. On the other hand, if $Q$ is the difference between the right and the left hand side in \eqref{eq:Hardy}, it has a generalized ground state $v(x)=|x|^\frac{2-n}{2}$ which is neither in $L^2(\R^n)$, nor in $\mathcal D^{1,2}(\R^n)$, nor in $L^2(\R^n;|x|^{-2})$. In this case, for $n>1$, $X^{**}$ is the subspace of $u\in W^{1,2}_{\mathrm{loc}}(\R^n)$  such that
$$
\int_{\R^n}|x|^{2-n}|\nabla(|x|^\frac{n-2}{2}u(x))|^2dx+\left(\int_{1<|x|<2}udx\right)^2<\infty,
$$
denoted in \cite{ky2} as $\mathcal D^{1,2}_{V}(\R^n)$ with $V(x)=-\left(\frac{n-2}{2}\right)^2\frac{1}{|x|^2}$.
\vskip5mm
\section*{Acknowledgments}
One of the authors (A.) thanks Department of Mathematics at Uppsala University for the hospitality. Theorem~\ref{thm:H2D} has been inspired by discussions of the other author (K.T.) with A.Tertikas and S.Filippas during his visit to Univsersity of Crete.

\end{document}